\documentclass[11pt]{article}
\usepackage{xspace,amssymb,amsmath,helvet,graphicx}
\begin{document}
\newcommand{\dee}{\,\mbox{d}}
\newcommand{\naive}{na\"{\i}ve }
\newcommand{\eg}{e.g.\xspace}
\newcommand{\ie}{i.e.\xspace}
\newcommand{\pdf}{pdf.\xspace}
\newcommand{\etc}{etc.\@\xspace}
\newcommand{\PhD}{Ph.D.\xspace}
\newcommand{\MSc}{M.Sc.\xspace}
\newcommand{\BA}{B.A.\xspace}
\newcommand{\MA}{M.A.\xspace}
\newcommand{\role}{r\^{o}le}
\newcommand{\signoff}{\hspace*{\fill} Rose Baker \today}
\newenvironment{entry}[1]%
{\begin{list}{}{\renewcommand{\makelabel}[1]{\textsf{##1:}\hfil}%
\settowidth{\labelwidth}{\textsf{#1:}}%
\setlength{\leftmargin}{\labelwidth}
\addtolength{\leftmargin}{\labelsep}
\setlength{\itemindent}{0pt}
}}%
{\end{list}}
\title{A useful family of fat-tailed distributions}
\author{Rose Baker\\School of Business\\University of Salford, UK\\
Email: rose.baker@cantab.net\\ORCID https://orcid.org/0000-0003-3555-3425}
\maketitle
\begin{abstract}
It is argued that there is a need for fat-tailed distributions that become thin in the extreme tail.
A 3-parameter distribution is introduced that visually resembles the t-distribution and interpolates between the normal distribution and the Cauchy distribution.
It is fat-tailed, but has all moments finite, and the moment-generating function exists.
It would be useful as an alternative to the t-distribution 
for a sensitivity analysis to check the robustness of results or for computations where finite moments are needed, such as in option-pricing.
It can be motivated probabilistically in at least two ways, either as the random thinning of a long-tailed distribution, or as random variation of the variance of a normal distribution.
Its properties are described, algorithms for random-number generation are provided,
and examples of its use in data-fitting given. Some related distributions are also discussed, including asymmetric and multivariate distributions.
\end{abstract}
\section*{Keywords}
Blurred truncation; Cauchy distribution; heavy tails; finite moments; poly-t distribution
\section{Introduction}
Heavy-tailed distributions that generalize the normal distribution abound, with the generalized hyperbolic distribution perhaps being the most all-inclusive, \eg Barndorff-Nielsen and Stelzer (2005).
However, in a finite world, heavy tails cannot continue for ever. A simple example is the Cauchy distribution (\eg Johnson {\em et al}, 1995), 
somewhat alarmingly exemplified as the distribution of bullet holes
made on a wall by a rotating gun. In the real world, the wall would be finite, giving sharp truncation to the Cauchy distribution. Even with an arbitrarily long wall, the bullets would not travel
arbitrarily far because of air resistance and gravity, giving censoring or what one might call `blurred truncation', a thinning of the tail. In general, one would expect blurred truncation of long-tailed distributions to be common in a finite world. 
Thus the distribution of the heights of athletes might have a long tail, but physical limits eventually impose a cut-off: the tallest man ever known of had a height of `only' 2.72 metres.

Thus ideally a fat-tailed distribution with an eventual cutoff is required, because
even a distribution with infinite thin tails has some probability for `impossible' values of the random variate. For example, concerning athlete heights, a thin-tailed distribution has a small probability
for the impossibly tall, or even for negative heights. However, estimating a cutoff parameter is difficult, and so it seems best to settle for
a distribution that is thin in the extreme tails and so has only a tiny probability for impossible occurrences.
The need for this more realistic modelling of long-tailed data is the motivation for this work.

There is a second reason for seeking fat-tailed distributions that are thin in the extreme tail.
Besides their physical implausibility, heavy-tailed distributions have some moments undefined or infinite, which causes problems for some analyses,
although usually the problem is only a minor one.
Thus long-tailed distributions model financial returns such as FTSE returns much better than the normal approximation used in deriving the famous Black-Scholes option pricing formula.
In fitting financial indices such as the FTSE, given successive values $S_{i-1}, S_i$, the `logged return' $X_i=\ln(S_1/S_{i-1})$ can be fitted to a t-distribution.
The actual return is $(S_i-S_{i-1})/S_{i-1}=\exp(X_i)-1$. Computing the expected return under some conditions then requires the existence of the means of exponentials of the random variable.
Hence Cassidy {\em et al} (2010) found that an integral needed for pricing European options (which was essentially the truncated moment-generating function) then diverged and
truncation of the t-distribution was needed to avoid an infinite result. 

In general, many workers have cited infinite moments as a problem and have introduced and studied truncated long-tailed distributions, \eg Nadarajah (2011). 
Apart from financial calculations of the kind described, where all moments must exist, in general the first four sample moments are useful for characterising distributions,
especially via measures such as skewness and kurtosis. They can also be used for estimating model parameters by the method of moments, or finding good starting values for methods of estimation such as maximum-likelihood. 
This last cannot be done if these moments are undefined or infinite
in the model. 

To construct a fat-tailed distribution with all moments existing, for the bullet-holes example given, assume a Raleigh distribution with survival function $S(x)=\exp(-\alpha x^2/2)$ for the distance $|x|$ along the wall a bullet can travel
in either direction.
The resulting pdf would be  $k\exp(-\alpha x^2/2)/(1+x^2)$, where the constant $k$ can be found analytically. 
It is here called the NC(1) distribution, and interpolates between a Cauchy distribution, as $\alpha\rightarrow 0$,
and a normal distribution as $\alpha\rightarrow\infty$.
The t-distribution also does so, although the t-distribution can have even longer tails than the Cauchy.  Here, however, all moments exist.
Further, the moment generating function also exists. 

This distribution is a special case of the `double-t' distribution. The poly-t distribution  pdf is a product of  $t$ pdfs, and is discussed by Dr\`{e}ze (1977), who gives an account of the origin of the name.
More recently Nadarajah (2009) has given the properties of the double-t distribution. However, the special cases given here are of value
because they are much more mathematically tractable than the general case, have fewer parameters, and have proved useful in fitting data.

The following sections describe the tractable types of double-t distribution. Finally some data fits are briefly discussed and conclusions drawn.
Many calculations have been checked via purpose-written  programs, and one such was used for the data fitting, using the NAG (Numerical Algorithms Group) library of routines.
\section{The double-t distribution and its special cases}
The special case of the poly-t distribution considered by Nadarajah (2009) has pdf
\[f(x)=C(1+x^2/a\alpha)^{-(1+a)/2}(1+x^2/b\beta)^{-(1+b)/2},\]
where the mean has been taken as zero.
This is the product of two t-pdfs with a common mean but different spreads and numbers of degrees of freedom.
Specialising to the case where $a \rightarrow \infty$ gives what one might call the N-t distribution, the general case of the distributions considered here.
This is a 4-parameter distribution where the long tails of the t-distribution can eventually become thin. 
Generalizing the Cauchy distribution to a t-distribution, a generalization of the NC(1) pdf  would have pdf
\begin{equation}f(x)=c\frac{\exp(-\alpha x^2/2)}{(1+x^2)^\gamma}\label{eq:gen}\end{equation}
for $\gamma > 0$, a t-type distribution with a Raleigh survival function for events to be accepted. 
We can write $\gamma=(\nu+1)/2$, where $\nu$ is the number of degrees of freedom of a t-distribution, and (\ref{eq:gen}) then generalizes the t-distribution.
However, it can show even more extreme tail behaviour when $0 < \gamma < 1/2$, before the eventual thin tail.
We can also take $\gamma < 0$, when the distribution becomes short-tailed. It is unimodal if $\gamma > \frac{\alpha}{2}$.
With $\gamma=-1/2$, the excess kurtosis is $\kappa \simeq -0.385$. This N-t distribution is then the short-tailed distribution of Tiku and Vaughan (1999).
With the affine transformation $X=(Y-\mu)/s$, the distribution of $Y$ has 4 parameters ($\mu, s, \alpha, \gamma$).

The constant $c(\alpha,\gamma)$ cannot be found analytically, and requires the evaluation of
an integral, but the distribution can still be used for data fitting, evaluating $c$ by numerical quadrature. 
This does not become a large computational burden even when $\mu$ is a function of covariates, because $c(\alpha,\gamma)$ 
only needs to be evaluated once per likelihood iteration as it is not a function of the mean $\mu$.
The parameters $\alpha, \gamma$ both control tail behaviour.

This distribution is fitted to data by maximum-likelihood later.  Some special cases are mathematically more tractable, the case where $\gamma=1/2$, given in appendix B, and where $\gamma$ is an integer. Calling these NC(n) distributions,
for pdfs that are the product of normal and Cauchy pdfs to the power $n$, the introductory example is an NC(1) distribution. The NC(2) distribution corresponds to a normal pdf
times a t[3] pdf, and so on. It was found as expected that the NC(1) distribution is needed for fitting very long-tailed data, but in practice the NC(2) distribution is arguably more 
useful, certainly for non-financial datasets. The NC(n) distributions are the main focus here, especially where $n$ is 1 or 2, but other cases are given in appendix A, and a multivariate distribution in appendix B.

\section{Properties of the NC(n) distributions}
\subsection{Probabilistic genesis and pdf}
The approach here can be thought of as taking a long-tailed distribution and censoring or thinning its tails according to the survival function of a second distribution,
to obtain a new distribution with some desired properties.
One probabilistic basis for the NC(1) distribution has already been alluded to: random numbers $X$ originate from a Cauchy distribution and are accepted with a probability $\exp(-\alpha X^2/2)$. 
Alternatively, a random number could originate from a normal distribution, and be accepted with
probability $1/(1+x^2)$.
Note that in the bullets on the wall example
the probability of acceptance varies with the line-of-flight distance $r$, which for a gun a distance $a$ from the wall is given by $r^2=x^2+a^2$.
Hence a Raleigh distribution for $r$ still gives a term $\exp(-\alpha x^2/2)$. 

Another probabilistic basis is to consider a normally-distributed random variable where the inverse standard deviation is a random variable with pdf $g(z)$, a Gaussian truncated below at $\alpha^{1/2}$, \ie $g(z)=\frac{\exp(-z^2/2)}{\sqrt{2\pi}\Phi(-\sqrt{\alpha})}$
for $z > \alpha^{1/2}$, where $\Phi$ is the normal distribution function.
Then the pdf is
\begin{equation}\frac{\int_{\alpha^{1/2}}^\infty z\exp(-(1+x^2)z^2/2)\dee z}{2\pi\Phi(-\sqrt{\alpha})}=k\frac{\exp(-\alpha x^2/2)}{1+x^2},\label{eq:pdf0}\end{equation}
where $k=\frac{\exp(-\alpha/2)}{2\pi\Phi(-\sqrt{\alpha})}$. 
The integral over the pdf (\ref{eq:pdf0}) is given in Gradshteyn and Ryzhik (2015) as result 3.466 (1), albeit in different notation.

It is the exclusion of high variances $> 1/\alpha$ that gives the blurred truncation to the Cauchy distribution.
With the generalized hyperbolic distribution, the generalizing distribution is the generalized inverse Gaussian distribution, so this does not produce blurred truncation.

The form (\ref{eq:pdf0}) is useful for derivations, but for practical use it is convenient to reparameterize the distribution, setting $\alpha=\beta/(1-\beta)$ and $X \rightarrow  X/(1-\beta)^{1/2}$, to obtain the pdf
\begin{equation} f(x)=c_1 \frac{\exp(-\beta x^2/2)}{1+(1-\beta)x^2},\label{eq:pdf}\end{equation}
where $c_1=\frac{(1-\beta)^{1/2}\exp(-\beta/2(1-\beta))}{2\pi\Phi(-{\beta/(1-\beta)}^{1/2})}$ and $0 < \beta \le 1$.
Then as $\beta\rightarrow 0$ the Cauchy distribution is obtained, and as $\beta\rightarrow 1$ the distribution becomes standard normal.
The random variable $Y$ will then be an affine transformation of $X$.
Without this reparameterization, the standard 1-parameter distribution would interpolate between Cauchy as $\alpha\rightarrow 0$ and a normal distribution of zero variance as $\alpha\rightarrow\infty$.

\begin{figure}[h]
\centering
\makebox{\includegraphics{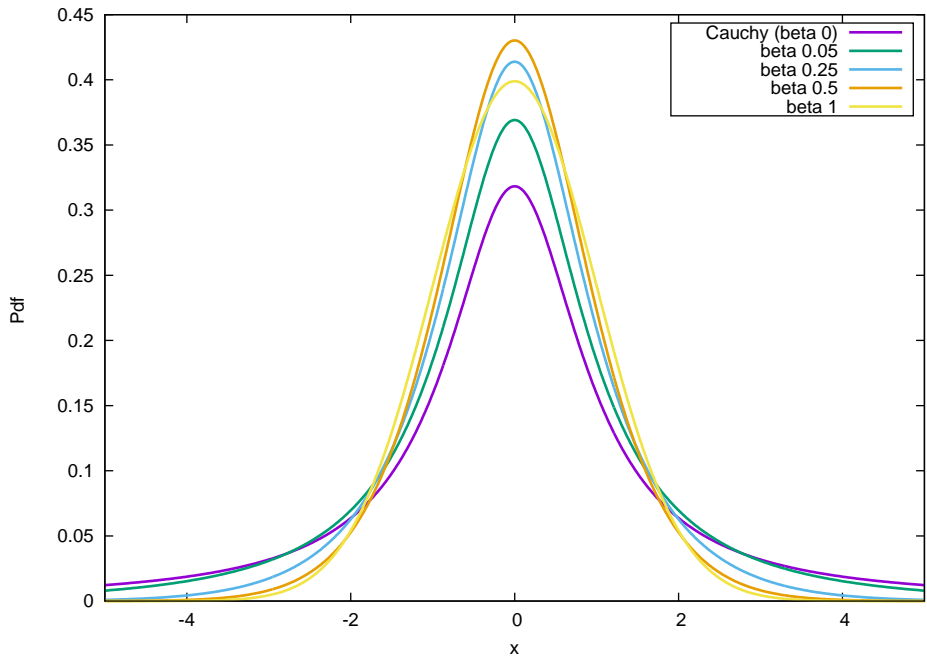}}
\caption{\label{figa}The pdf of the NC(1) distribution for different values of the tail-weight parameter $\beta=\alpha/(1+\alpha)$.}
\end{figure}

\begin{figure}[h]
\centering
\makebox{\includegraphics{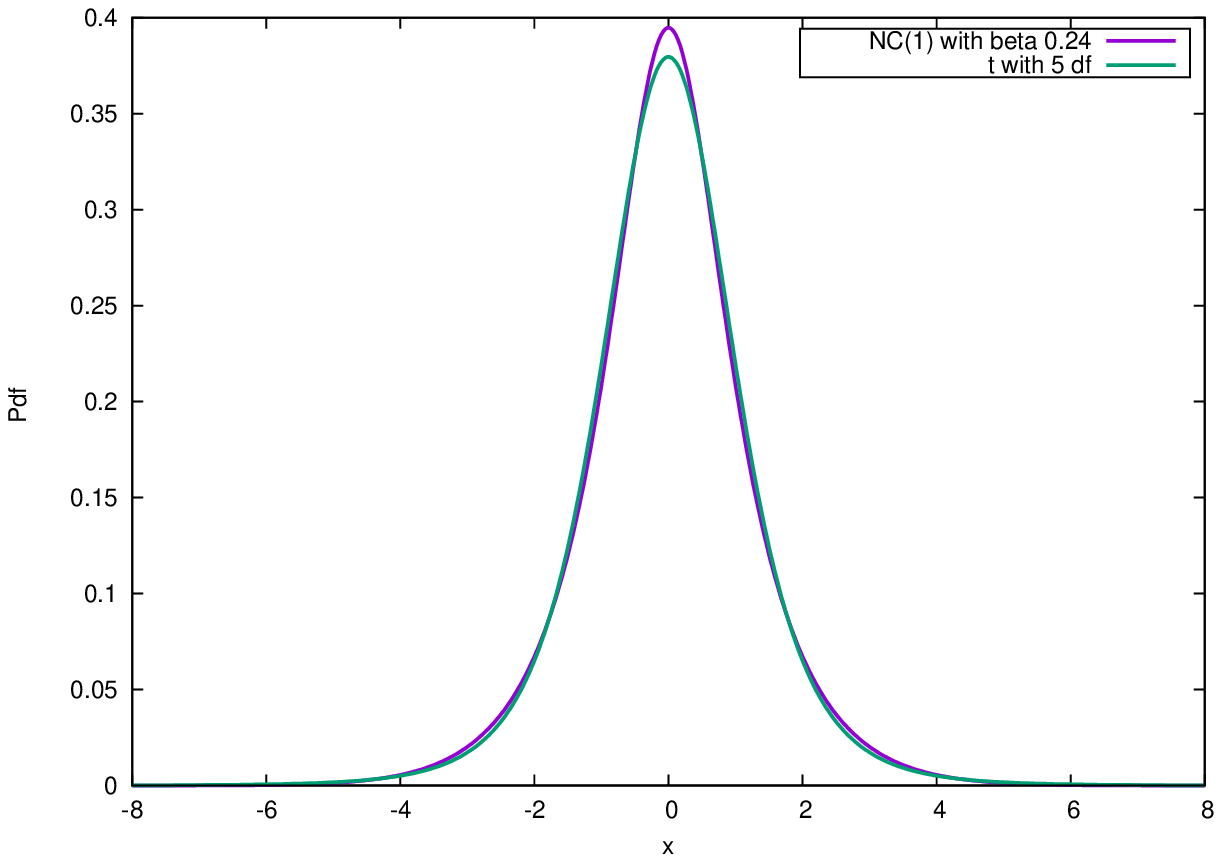}}
\caption{\label{figb}The pdf of the best-matching NC(1) distribution with the t-distribution with 5 degrees of freedom, with match measured by Hellinger distance.}
\end{figure}

In general, $\beta$ here is analogous to the number of degrees of freedom $\nu$ of the t-distribution.
Figure \ref{figa} shows the NC(1) pdf for several values of $\beta$, 
and figure \ref{figb} shows the  t-distribution with $\nu=5$ degrees of freedom with the closest NC(1) distribution, as measured by Hellinger distance. 

This was selected for illustration as the distributions are identical with $\nu=1$ and $\nu$ infinite, but it was not evident how close they would be in between.
Except in the extreme tails, the distributions look very similar. The NC(1) distribution has a slightly higher peak, because there is less probability in the extreme tails.
It has slightly thinner shoulders, and is fatter in the tail before the tail finally thins. The excess kurtosis is 6 for the t-distribution, and 1.50 for the NC(1) distribution,
showing that although the distributions look similar, the NC(1) tail is eventually thin.

The `standard' 1-parameter form for $X$ is described here, with $X = (Y-\mu)/s$ for the 3-parameter distribution of $Y$ with arbitrary mean and variance.
\subsection{Moments}
\begin{figure}[h]
\centering
\makebox{\includegraphics{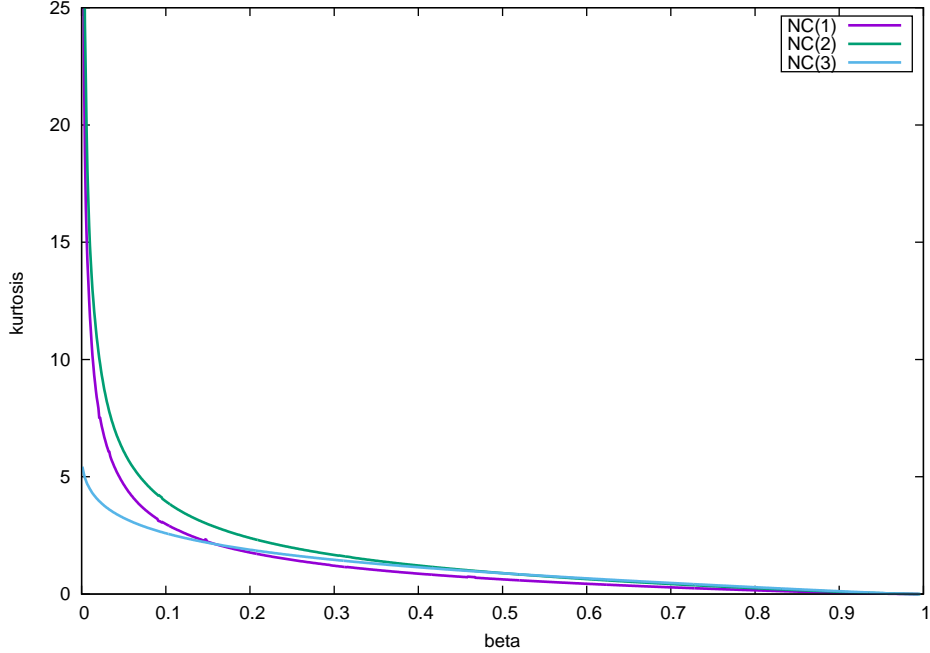}}
\caption{\label{figd}Excess kurtosis versus $\beta$ for NC(1), NC(2) and NC(3) distributions.}
\end{figure}

Moments exist for $\beta > 0$.
The odd moments are all zero by symmetry. Evaluation of the integrals for second and fourth moments $\sigma^2=\text{E}(X^2)$ and $\text{E}(X^4)$
can be done given (\ref{eq:pdf}) and the normal integral. Then
\[\text{E}(X^2)=(1-\beta)^{-1}\{\frac{c_1\sqrt{2\pi}}{\beta}-1\},\]
\[\text{E}(X^4)=(1-\beta)^{-2}\{c_1\frac{\sqrt{2\pi/\beta}(1-2\beta)}{\beta}+1\}.\]
The excess kurtosis goes from infinity at $\beta=0$ to zero at $\beta=1$,
as shown in figure \ref{figd}.

The moment generating function $M(t)=\text{E}\{\exp(tX)\}$ is the integral
\[M(t)=\frac{1}{\sqrt{2\pi}\Phi(-\sqrt{\beta/(1-\beta)})}\int_{\{\beta/(1-\beta)\}^{1/2}}^\infty \exp\{\frac{t^2}{2(1-\beta)z^2}-z^2/2\}\dee z,\]
from which moments may be derived by expanding the first term in the exponent. The mgf exists for $\beta > 0$.

The characteristic function $C(t)$ is the same but with $t^2$ replaced by $-t^2$. The integral is related to the leaky aquifer function, and can be written as
\[C(t)=\frac{1}{2\sqrt{\pi}\Phi(-\sqrt{\alpha})}\int_\alpha^\infty y^{-1/2}\exp\{-y-t^2/(4(1-\beta)y)\}\dee y\]
\[=\frac{1}{2\sqrt{\pi}\Phi(-\sqrt{\alpha})}\Gamma\{1/2,\alpha;t^2/(4(1-\beta))\},\]
where $\Gamma$ denotes the generalized incomplete gamma function. When the first argument is half-integer, as here, this can be written in closed form, in terms of special functions (Chaudhry {\em et al}, 1996).
It would however be easier and probably much quicker to evaluate the integral numerically.

\subsection{Distribution function}
This cannot be written in closed-form, but must be obtained by integrating (\ref{eq:pdf}).
A useful result follows from reversing the order of the integrations over $x$ and $z$.
Then we have that the distribution function $F(x)$ is given by
\[F(x)=\frac{1}{\sqrt{2\pi}\Phi(-\sqrt{\alpha})}\int_{\alpha^{1/2}}^\infty \Phi(\sqrt{1-\beta}z x)\exp(-z^2/2)\dee z.\]
This form is arguably less convenient for numerical quadrature, but when $X=(1-\beta)^{-1/2}$ the integral can be done, to give the results
\[F((1-\beta)^{-1/2})=1-\Phi(-\sqrt{\alpha})/2,\]
\[F(-(1-\beta)^{-1/2})=\Phi(-\sqrt{\alpha})/2.\]
These results assist numerical integration, which can now be done over a smaller range.


\subsection{Random numbers}
These can be generated in several ways using the rejection method. A simple and fairly efficient method is to generate a random variable $X$ from a normal distribution,
to generate $Z$ from a truncated normal distribution, and then to set $X \rightarrow \frac{X}{(\sqrt{1-\beta})Z}$.
Robert (1995) used an exponential majorizing distribution in a rejection method to generate random numbers from the truncated normal distribution.
For completeness, the full method including his algorithm is given here.
\begin{enumerate}
\item Set $\lambda=(\sqrt{\alpha}+\sqrt{4+\alpha})/2$;
\item generate $U, V$,  uniformly-distributed random numbers;
\item generate the shifted exponential random variable $Z=\sqrt{\alpha}-\ln(U)/\lambda$;
\item compute $\rho=\exp(-(Z-\lambda)^2/2)$;
\item if $V > \rho$ restart from step 2;
\item generate $X$, a normally-distributed random number;
\item return the random number as $\frac{X}{(\sqrt{1-\beta})Z}$.
\end{enumerate}
\subsection{Estimation}
With $Y$ observed, where $(y-\mu)/s=x$, the pdf of $Y$ is
\begin{equation} f(y)=(c_1/s) \frac{\exp(-\beta (y-\mu)^2/2s^2)}{1+(1-\beta)(y-\mu)^2/s^2}.\label{eq:pdfy}\end{equation}
Write $D_i=1+(1-\beta)\{(y_i-\mu)/s\}^2$,
then on making $n$ observations $y_1\cdots y_n$, the log-likelihood is
\begin{equation}\ell=n\ln c_1-n\ln s -(\beta/2s^2)\sum_{i=1}^n (y_i-\mu)^2-\sum_{i=1}^n \ln D_i.\end{equation}
Then with
\begin{align*}
\ln c_1&=(1/2)\ln(1-\beta)-\frac{\beta}{2(1-\beta)}-\ln(2\pi)-\ln\Phi\{-\sqrt{\frac{\beta}{1-\beta}}\},\\
\dee \ln c_1/\dee \beta&=-(1/2)\{\frac{1}{1-\beta}+\frac{1}{(1-\beta)^2}\}+\frac{\exp\{-\frac{\beta}{(2(1-\beta)}\}}{2(2\pi)^{1/2}\Phi\{-\sqrt{\frac{\beta}{(1-\beta)}}\}\beta^{1/2}(1-\beta)^{3/2}},\\
\dee^2 \ln c_1/\dee\beta^2&=-\frac{1}{2(1-\beta)^2}-\frac{1}{(1-\beta)^3}\\
&-\frac{\exp\{-\frac{\beta}{(2(1-\beta)}\}}{4(2\pi)^{1/2}\Phi\{-\sqrt{\frac{\beta}{1-\beta}}\}}\{\beta^{-1/2}(1-\beta^{-7/2})+\beta^{-3/2}(1-\beta)^{-3/2}+3\beta^{-1/2}(1-\beta)^{-5/2}\}\\
&+\frac{\exp(-\beta/(1-\beta))\beta^{-1}(1-\beta)^{-3}}{8\pi\Phi\{-\sqrt{\frac{\beta}{(1-\beta)}}\}^2}.
\end{align*}
the likelihood derivatives are
\begin{align*}
\partial\ell/\partial\mu&=(\beta/s^2)\sum_{i=1}^n (y_i-\mu)+\frac{2(1-\beta)}{s^2}\sum_{i=1}^n (y_i-\mu)/D_i,\\
\partial\ell/\partial s&=-n/s+(\beta/s^3)\sum_{i=1}^n (y_i-\mu)^2+2\{(1-\beta)/s^3\}\sum_{i=1}^n (y_i-\mu)^2/D_i,\\
\partial\ell/\partial\beta&=n\dee\ln(c_1)/\dee \beta-(1/2s^2)\sum_{i=1}^n (y_i-\mu)^2+(1/s^2)\sum_{i=1}^n (y_i-\mu)^2/D_i,
\end{align*}
and the second derivatives are
\begin{align*}
\partial^2\ell/\partial\mu^2&=-n\beta/s^2+\{2(1-\beta)/s^2\}\{2(1-\beta)/s^2\sum_{i=1}^n (y_i-\mu)^2/D_i^2-\sum_{i=1}^n 1/D_i\},\\
\partial^2\ell/\partial\mu\partial s&=-2(\beta/s^3)\sum_{i=1}^n (y_i-\mu)-4\{(1-\beta)/s^3\}\sum_{i=1}^n (y_i-\mu)/D_i\\
&+4\{(1-\beta)^2/s^5\}\sum_{i=1}^n (y_i-\mu)^3/D_i^2,\\
\partial^2\ell/\partial\mu\partial\beta&=(1/s^2)\sum_{i=1}^n (y_i-\mu)-(2/s^2)\sum_{i=1}^n (y_i-\mu)/D_i+2\{(1-\beta)/s^4\}\sum_{i=1}^n (y_i-\mu)^3/D_i^2,\\
\partial^2\ell/\partial s^2&=n/s^2-3(\beta/s^4)\sum_{i=1}^n (y_i-\mu)^2-6\{(1-\beta)/s^4\}\sum_{i=1}^n (y_i-\mu)^2/D_i\\
&+4\{(1-\beta)^2/s^6\}\sum_{i=1}^n (y_i-\mu)^4/D_i^2,\\
\partial^2\ell/\partial s\partial\beta&=(1/s^3)\sum_{i=1}^n (y_i-\mu)^2-(2/s^3)\sum_{i=1}^n (y_i-\mu)^2/D_i+\{2(1-\beta)^2/s^5\}\sum_{i=1}^n (y_i-\mu)^4/D_i^2,\\
\partial^2\ell/\partial\beta^2&=n\dee^2\ln(c_1)/\dee \beta^2+(1/s^4)\sum_{i=1}^n (y_i-\mu)^4/D_i^2.
\end{align*}

This is an instance where the NC(1) distribution is more tractable than the t-distribution, because differentiating the corresponding constant for the t-distribution
would require the computation of the psi (digamma) function.
\subsection{NC(n) distributions}
For distributions of the form
\begin{equation}f(x)=c_n\frac{\exp(-\alpha x^2/2)}{(1+x^2)^n},\label{eq:pdfn}\end{equation}
for integer $n$, the constant $c_n$ can be found analytically. Such distributions can be useful in fitting \eg financial data, and so the derivation of the constant $c_n$ is outlined.
The generalizing function $g(z)$ is taken as 
\[g(z)=d_n (z^2-\alpha)^{n-1}\exp(-z^2/2)\]
for $z \ge \sqrt{\alpha}$. This is 
chosen to give an integral of form (\ref{eq:pdfn}).
Then 
\begin{align}f(x)&=d_n\int_{\alpha^{1/2}}^\infty\frac{z\exp(-(1+x^2)z^2/2)(z^2-\alpha)^{n-1}\dee z}{\sqrt{2\pi}}\nonumber \\
&=(\frac{d_n2^{n-1}(n-1)!\exp(-\alpha^2/2)}{\sqrt{2\pi}})\frac{\exp(-\alpha x^2/2)}{(1+x^2)^n},\label{eq:fn}\end{align}
where the last result follows on integrating by parts $n$ times.

It remains to find the constant $d_n$, where $d_n^{-1}=I_n=\int_{\sqrt{\alpha}}^\infty (z^2-\alpha)^{n-1}\exp(-z^2/2)\dee z$.
Using integration by parts yet again on $g(z)$ , we obtain
\[I_1=\sqrt{2\pi}\Phi(-\sqrt{\alpha}),\]
\[I_2=\sqrt{\alpha}\exp(-\alpha/2)+(1-\alpha)\sqrt{2\pi}\Phi(-\sqrt{\alpha}).\]
Then for $n > 2$ integration by parts gives the recurrence relation
\[I_n=(2n-3-\alpha)I_{n-1}+2\alpha(n-2)I_{n-2}.\]
Hence it is possible to find the constants $c_n$ from $d_n=I_n^{-1}$ and (\ref{eq:fn}) and fit these distributions to data.
The corresponding number of degrees of freedom is $\nu=2n-1$. Again $X \rightarrow X/\sqrt{1-\beta}$ for practical use.

Random numbers can be generated from the pdf
\[f(x)=\frac{\sqrt{1-\beta}c_n \exp(-\beta x^2/2)}{(1+(1-\beta)x^2)^n}\]
using the rejection method, by generating a random number from the $t[2n-1]$ distribution, rescaling it, and accepting it with probability $\exp(-\beta x^2/2)$.
In detail:
\begin{enumerate}
\item Generate a random number $Y$ from the $t[2n-1]$ distribution, \eg using the method of Bailey (1994);
\item Rescale to $X=Y/\sqrt{(2n-1)(1-\beta)}$;
\item Generate a uniform random number $U$ and accept if $U < \exp(-\beta X^2/2)$, else return to step 1.
\end{enumerate}
This algorithm can also be used for the NC(1) distribution.
Timing tests showed that for the NC(1) distribution, the random-number algorithm given specifically for that distribution took only 54.2\% as much time
as the general algorithm here. This was based on a range of values of $\beta$ from 0.01 to 0.99.

\subsection{The NC(2) distribution}
This distribution cannot attain the Cauchy-like tails of the  NC(1) distribution, but has proved useful in fitting the datasets described later; even financial datasets
are not as long-tailed as the Cauchy.

The pdf is 
\[f(x)=\frac{\sqrt{1-\beta}c_2 \exp(-\beta x^2/2)}{(1+(1-\beta)x^2)^2},\]
with $c_2$ defined as in the previous section, \ie
\[c_2=\frac{2}{\sqrt{2\pi\alpha}+2\pi(1-\alpha)\exp(\alpha/2)\Phi(-\sqrt{\alpha})}.\]

The first two even moments can be easily derived and are:
\[\text{E}(X^2)=(1-\beta)^{-1}(c_2/c_1-1),\]
\[\text{E}(X^4)=(1-\beta)^{-2}(1-2c_2/c_1)+(1-\beta)^{-3/2}\sqrt{2\pi/\beta}c_2.\]
The excess kurtosis as a function of $\beta$ is shown in figure \ref{figd}. Note that as $\beta\rightarrow 0$ the kurtosis is infinite for NC(1) and NC(2) distributions.

Random numbers can be generated as previously described for the NC(n) distribution.

\section{Fits to data}
The fitting was done using a purpose-written Fortran program and using the NAG library function minimizers and quadrature routines. Derivatives were not required by the gradient-based function minimisers used, and identical results were obtained with 
different minimisers. The model parameters were the mean $\mu$, scale factor $s$, and tail parameters $\beta$ and/or $\nu$.
Iteration converged with no problem using the sample mean and standard deviation as starting values for $\mu, s$ and starting estimates of $\beta, \nu$ at middle-of the range values, \eg
$\beta=0.2, \nu=5$. However, random restarts were made by perturbing the parameter values from their fitted values and re-minimising, to check that the global maximum of the log-likelihood had been attained.

Several publicly-available datasets were fitted. There is little point showing histograms with the fitted model, as the fits are visually very similar to that of the t-distribution,
as can be seen from figure \ref{figb}.
First, the heights in centimetres of 100 female athletes, from Cook and Weisberg (1994).
Next, the 1080 monthly-average heights of the Rio Negro river at Manaus (Sternberg, 1987), the 10939 logged returns on the S \& P 500 index to
20 may 2021, the 9013 logged returns on the FTSE-100 to 8 April 2021, and the 13863 logged returns on the Nikkei to 20 May 2021.
The first two datasets are physical measurements, and the last three are financial.
It is known that day-to-day returns are not independent, so a better fit to data can be obtained by modelling this dependence, but this was ignored here
for simplicity.

Results  from fitting the $t$, NC(1), NC(2) and N-t distributions are shown in table \ref{tab1}. 
\begin{table}[h]
\begin{tabular}{|l|c|c|c|c|} \hline
Dataset&Method &$-\ell$&$\beta$&$\nu$ \\ \hline
Athletes  & $t$ &349.36&-& 4.24 (2.080) \\ \hline
(100)&NC(1)&348.77&0.18 (.156) & - \\ \hline
&NC(2) & 349.09&0.175 (.275) & - \\ \hline
&N-t & 348.76 & 0.166 (.171) & 0.71 (1.747) \\ \hline \hline
River height  & $t$ &1974.45 & - & 6.43 (1.234) \\ \hline
(1080)&NC(1) & 1975.46&0.323 (.097) & - \\ \hline
&NC(2)&1974.16&0.337 (.137) & - \\ \hline
&N-t&1974.10&0.309 (.236)&1.34 (.918) \\ \hline \hline
S \& P 500  & $t$&15143.32&-&2.95 (.097) \\ \hline
(10939)&NC(1)&15220.70&0.050 (.004) & - \\ \hline
&NOD (2) & 15142.90&0.0038 (.004)&-\\ \hline
&N-t&15141.37&0.0124 (.007) & 2.73 (.156) \\ \hline\hline
FTSE-100  &$t$&12778.91&-&3.51 (.140) \\ \hline
(9013)&NC(1)&12857.36&0.088 (.007) & - \\ \hline
&NC(2)&12780.41&0.029 (.010) & - \\ \hline
&N-t&12778.65&0.0087 (.014) & 3.40 (.211) \\ \hline\hline
Nikkei  & $t$& 21410.39&-&2.88 (.086) \\ \hline
(13863)&NC(1) &21407.52&0.058 (.004)&-\\ \hline
&NC(2)&21407.82&0.0122 (.006) & - \\ \hline
&N-t&21386.32&0.0504 (.007) & 1.88 (.152) \\ \hline

\end{tabular}
\caption{Minus the log-likelihood and parameters $\beta, \nu$ with standard error for the 5 datasets studied.
The models are the t-distribution, Normal-Cauchy, NOD of degree two and normal-t distribution\label{tab1}.
}
\end{table}
The parameters $\mu, s$ are of no interest for our purpose and are not shown. All distributions except the N-t have 3 parameters and so log-likelihoods can be compared. The N-t distribution has 4 parameters.

Overall, the NC(1) distribution fits comparably with the t-distribution, sometimes better and sometimes worse (better for the height of athletes data and the Nikkei).
The NC(2) distribution fits better than the t-distribution in 4 out of the 5 cases. The N-t distribution generalizes the t-distribution and so must always fit better in the sense that the log-likelihood will be higher.
The improvement is large for the Nikkei, where the log-likelihood increases by over 24. The large size of the dataset enables small differences in fit to be detected. Because $\nu$ is somewhat lower for the N-t fit, it seems that the Nikkei
returns favour a very long-tailed distribution, but one where the tails do eventually attenuate. 
It is clear that the N-t distribution gives a significantly better fit than the t-distribution, so that $X^2[1]=48.14$.

The aim here was to show `noninferiority', \ie that the NC(n) distributions could data comparably to the t-distribution, sometimes slightly better and sometimes slightly worse.
This is the case here. 
The tentative conclusion is the NC(2) distribution may be more useful in practice than the NC(1)
distribution. Although it does not allow such extreme tail behaviour, it is better at modelling the less extreme tail behaviour that usually occurs.
However, the NC(1) distribution allows extreme tail behaviour and will always be adequate.

It is interesting to examine the excess kurtosis of the three financial indices, where with the t-distribution, the predicted kurtosis would be infinite. Data were divided into 10 tranches, and randomized, to obtain sample estimates $\hat{\kappa}$,
while predictions $\kappa_{\text{pred}}$ were made from the fitted N-t model.
The results were: S \& P 500: $\hat{\kappa}=19.42 \pm 9.04, \kappa_{\text{pred}}=14.87$, FTSE-100 $\hat{\kappa}=9.32 \pm 1.60, \kappa_{\text{pred}}=12.68$,
Nikkei $\hat{\kappa}=8.055 \pm 2.11, \kappa_{\text{pred}}=6.33$. This shows that the predicted kurtoses are in fair agreement with the sample kurtoses, confirming that the
new models give a better description of the data than the t-distribution.

\section{Conclusions}
A class of fat-tailed distributions that are normal in the extreme tails has been described, with a number of generalizations and related distributions.
Figure \ref{figb} shows the NC(1) distribution to be very similar
in shape to the t-distribution, except in the extreme tails, and it has been demonstrated to fit datasets very comparably to the t-distribution.
The data-fitting led to the recognition that the NC(2) distribution also looks promising for practical work, as it fitted 4 out of 5 datasets better than the t-distribution.

Compared to the t-distribution, the computation of the pdf is of comparable difficulty, both requiring special functions for the normalization constant.
Moments are simply computed, just being somewhat messier than those of the t-distribution. The moment generating function also exists, and can be found as an integral.
Random numbers can be easily generated, which is vital for many computations, \eg Markov-chain Monte-Carlo.
The distribution function must be found by numerical quadrature, whereas for the t-distribution it is available as a special function, the incomplete beta function.

One might ask, given that we have the t-distribution, why we should be interested in a similar distribution, where the computations are slightly more complex.
There are at least three reasons. Firstly, from the arguments given in the introduction, we expect that these distributions will model real-world
behaviour better in the extreme tails, which could be useful for predicting the probability of extreme events.

Secondly, for sensitivity analysis: to check robustness of results to model assumptions we can repeat our analysis with the NC(1) distribution replacing the t-distribution,
and we should still reach similar conclusions.

Thirdly, an advantage of the NC(1) distribution as a replacement for the t-distribution is that all the moments and the moment-generating function exist. This has two benefits. One is that the method of moments can be used to
find good starting values for a maximum-likelihood fit to the data, \eg by using figure \ref{figd}. The main benefit however is for cases where the mean of the exponent of the random variable must be computed.
This happens with long-tailed data such as financial return data, where logarithms have been taken and the distribution is still long-tailed. The mean of the exponent then gives the mean of the returns themselves.
This situation also occurs where for example BMI (body mass index) or other medical statistics might be fitted by taking logarithms. It will be the mean and variance of BMI itself, not its logarithm, that are ultimately required.
 
Future work in this area could be gaining experience on fitting the distribution (\ref{eq:pdf}) in different application areas, and investigating further related distributions.
Random number generation is a practically useful topic, and the generation schemes given here could doubtless be honed.
\section*{Appendix A: a multivariate distribution}
\subsection{Multivariate distributions}
Writing $A=({\bf x}-{\boldsymbol \mu})^T{\bf V}^{-1}({\bf x}-{\boldsymbol \mu})$, 
where ${\boldsymbol \mu}$ is the mean and ${\bf V}$ the covariance matrix of the normal distribution, and with all variances scaled by $z^2$ as before,
the p-dimensional pdf becomes
\[f({\bf x})=\frac{1}{(2\pi)^{p/2}|{\bf V}|^{1/2}}\int_{\alpha^{1/2}}^\infty z^p \exp(-(1+A)/2)g(z)\dee z.\]
To obtain the desired pdf, take $g(z)=c_z \exp(-z^2/2)/z^{p-1}$.
Then for $p=2$, the bivariate case,  $g(z)=c_z \exp(-z^2/2)/z$, and
$c_z=2/E_1(\alpha/2)$. The required pdf is therefore
\[f({\bf x})=\{\frac{1}{\pi |{\bf V}|^{1/2} E_1(\alpha/2)}\}\frac{\exp(-\alpha (1+A)/2)}{1+A}.\]
For the trivariate case, $g(z)=c_z \exp(-z^2/2)/z^2$ and the constant is $c=c_z (2\pi)^{-3/2}|{\bf V}|^{-1/2}$, where 
\[c_z=\frac{1}{\exp(-\alpha/2)/\sqrt{\alpha}-\sqrt{2\pi}\Phi(-\sqrt{\alpha})}.\]

In general, taking ${\bf V}/z^2$ as the covariance matrix of a normal distribution, the covariance matrix $\bf \Sigma$ is
\[{\bf \Sigma}={\bf V}c_z\int_{\alpha^{1/2}}^\infty \frac{\exp(-z^2/2)\dee z}{z^{p+1}}.\]
By parts,
\[{\bf \Sigma}={\bf V}\{\frac{c_zp\exp(-\alpha/2)}{\alpha^{p/2}}-1\}.\]

\section*{Appendix B: other cases}
\subsection{A blurred t-distribution with zero degrees of freedom}
When $\gamma=1/2$, a different  distribution is obtained.
The modified Bessel function of the second kind can be written as 
\[K_\nu(z)=\frac{\sqrt{\pi}(z/2)^\nu}{\Gamma(\nu+1/2)}\int_0^\infty \exp(-z \cosh t)\sinh^{2\nu}(t) \dee t.\]
Setting $\nu=0, x=\sinh (t/2), z=\alpha/4$ we obtain the pdf
\[f(x)=\{\frac{\exp(-\alpha/4)}{K_0(\alpha/4)}\}\frac{\exp(-\alpha x^2/2)}{\sqrt{1+x^2}}.\]
This gave similar fits to (\ref{eq:pdf0}) for the datasets fitted here, but is less tractable.
Moments can be computed in terms of the Bessel functions $K_1, K_2$ \etc
\subsection{Asymmetric distributions}
To introduce asymmetry to the NC(1) distribution, a simple method is to use a 2-piece distribution:

\[f(x)= \left\{\begin{array}{ll}
\frac{2k}{s_1+s_2} \frac{\exp(-\beta (x-\mu)^2/2s_1^2)}{1+(1-\beta)(x-\mu)^2/s_1^2} & \mbox{if $X \le \mu$} \\
\frac{2k}{s_1+s_2} \frac{\exp(-\beta (x-\mu)^2/2s_2^2)}{1+(1-\beta)(x-\mu)^2/s_2^2} & \mbox{if $X > \mu,$ }
\end{array}
\right. \]
where $s_1, s_2$ are scale factors. The pdf of this type of distribution and its first derivative are continuous at $X=\mu$.
The mean can be shown to be
\[\text{E}(X)=\mu+k(s_2-s_1)\exp(\frac{\beta}{2(1-\beta)})E_1( \frac{\beta}{2(1-\beta)})/2(1-\beta),\]
where $E_1$ is the exponential integral $E_1(z)=\int_z^\infty \frac{\exp(-t)\dee t}{t}$.

For $\beta > 0$, even moments about $\mu$ are simply derived, \eg the second moment is $(s_1^3+s_2^3)/(s_1+s_2)$ times the second moment of (\ref{eq:pdf}).
Thus the moments can be derived in terms of special functions, albeit messily. The probability that $X \le \mu$ is $s_1/(s_1+s_2)$, so random numbers can be generated
by finding which half of the range  $X$ lies in, then proceeding as for the NC(1) distribution.

\subsection{Survival distributions}
The two distributions used in this approach can of course be survival distributions, defined for $X \ge 0$.
Tebbens {\em et al} (2001) used a distribution $Cr^{-\alpha}\exp(-ar)$ to model the tail behaviour of the volume of seamounts.

Let $f(x)=c\frac{\exp(-\beta x)}{\{1+(1-\beta)x\}^\gamma}$. This distribution interpolates between the Lomax distribution when $\beta=0$ and the exponential when $\beta=1$,
and has three parameters, on letting $X$ be rescaled. The constant
\[c=\frac{\beta\exp(-\beta/(1-\beta))(\beta/(1-\beta))^{-\gamma}}{\Gamma(1-\gamma,\beta/(1-\beta))},\]
where the incomplete gamma function is $\Gamma(\alpha,x)=\int_x^\infty z^{\alpha-1}\exp(-z)\dee z$.

\section*{Acknowledgements}
The author would like to thank Professors Philip Scarf and Ian McHale for helpful comments.


\begin{thebibliography}{99}
\bibitem{bailey}Bailey, R. W. (1994). Polar generation of random variates with the t-distribution, Mathematics of Computation, {\bf 62} (206), 779-781.
\bibitem{barn}Barndorff-Nielsen O. E. and Stelzer, R. (2005). Absolute moments of generalized hyperbolic distributions and approximate scaling of normal inverse Gaussian L\'{e}vy processes, 
Scandinavian Journal of Statistics {\bf 32} (4), 617-637.
\bibitem{cass}Cassidy, D.T., Hamp, M.J. and Ouyed, R. (2010). Pricing European
options with a log Student's t-distribution: a Gosset formula. Physica A, {\bf 389}, 5736-5748.
\bibitem{chaud}Chaudhry, M. A., Temme, N. M. and Veling, E. J. M. (1996). Asymptotics and closed form of a generalized incomplete gamma function, Journal of Computational and Applied Mathematics,
{\bf 67}, 371-379.
\bibitem{cook}Cook, R. D. and Weisberg, S. (1994). An Introduction to Regression Graphics. Wiley, New York.
\bibitem{dreze}Dr\`{e}ze, J. H. (1977). Bayesian regression analysis using poly-t densities, Journal of Econometrics {\bf 6}, 329-345.
\bibitem{grad} Gradshteyn I. S. and Ryzhik I. M. (2015). Table of integrals, series, and products, 8th ed., Academic Press, Waltham.
\bibitem{kotz} Johnson, N. L., Kotz, S. and Balakrishnan, N. (1995), Continuous univariate distributions vol. 1, Wiley, New York.
\bibitem{nada0}Nadarajah, S. (2009). The product $t$ density distribution arising from the product of two student's $t$ PDFs, Statistical papers {\bf 50} 605-615.
\bibitem{nada}Nadarajah, S. (2011). Making the Cauchy work, Brazilian Journal of Probability and Statistics,{\bf 25} (1),99-120.
\bibitem{rob}Robert, C. P. Simulation of truncated normal variables (1995). Statistics and Computing {\bf 5}, 121-125.
\bibitem{stern}Sternberg, H. O'R. (1987) Aggravation of floods in the Amazon river as a consequence of deforestation? Geografiska Annaler, {\bf 69}A, 201-219. 
\bibitem{tebbens}Tebbens, S. F., Burroughs, S. M., Barton, C. C. and Naar, D. F. (2001). Statistical self-similarity of hotspot seamount volumes modeled as self-similar criticality,
Geophysical Research Letters, {\bf 28} (14), 2711-2714.
\bibitem{tiku}Tiku, M. L. and Vaughan, D. C. (1999). A family of short-tailed symmetric distributions. Technical report, McMaster University, Canada.
\end{thebibliography}
\end{document}